\input amstex
\documentstyle{amsppt}
\NoRunningHeads
\define\C#1{\Cal #1}

\define\norm#1#2{\|{#2}\|_{#1}}
\define\Int{\int\limits}
\define\Hint#1{\Int_{\C_-} #1\mu(dx)}

\define\Sum{\sum\limits}

\define\id{\operatorname{id}}
\define\Id{\operatorname{Id}}
\define\Eint#1#2#3{\Int_{#1}^{#2}#3\widehat{dw(\tau)}}
\define\pair#1#2#3{\langle #2,\, #3\rangle_{#1}}
\define\CL{\C{L}}
\define\F{\frak F}
\TagsOnRight
\topmatter
\title
Stochastic Equations and Evolution Families in the space of Formal Mappings
\endtitle
\author
I.Ya. Spectorsky \\ National Technical University of Ukraine (KPI),
Kiev, Ukraine
\endauthor
\thanks
This work was supported, in part, by the International
Soros Science Education Program  (ISSEP) through grant  N PSU051117.
\endthanks
\endtopmatter

\document

\head 1. Introduction \endhead

Let's give some consideration leading to the notion of formal mapping.

Let's consider a stochastic equation
$$
y(t)=y(0)+\Int_0^t a(\tau)(y(\tau))d\tau
+\Int_0^t b(\tau)(y(\tau))dw(\tau),\qquad 0\leqslant t\leqslant T
$$
in Hilbert space $Y$ (all Hilbert spaces are supposed to be real and
separable).  Here $w$ is a Wiener process, associated with
 canonical triple $H_+\subset H_0\subset
H_-$,  with Hilbert-Schmidt embeddings, $a$ and $b$ are
continuous mappings from $[0,T]\times Y$ into spaces $Y$ and
$\CL_2(Y)$ respectively, $y$ is unknown random process taking its
values in space $Y$, $y(0)$ is nonrandom initial condition.

Let coefficients $a$ and $b$ be analytical functions with respect to $y
\in Y$ under fixed $t$, and $a(t)(0)=0$, $b(t)(0)=0$. I.e., for all
$y\in Y$ the following expansions into power series holds:
$$
\gather
a(t,y)=\Sum_{k\geqslant1} a_k(t)(y,y,\dots,y),\quad
b(t,y)=\Sum_{k\geqslant1} b_k(t)(y,y,\dots,y).
\endgather
$$
Here $a_k(t)$ and $b_k(t)$ are $k$-linear continuous operators from $Y$ to
$Y$ and from $Y$ to $\CL_2(Y)$ respectively.

Let this equation have the unique solution $y(t)=S(t)(y)$,
and the following expansion holds: $S(t)(y)=\Sum_{k\geqslant1}
S_k(t)(y,y,\dots,y)$. In this case, function
$a(t)\circ S(t)(y)$ and $b(t)\circ S(t)(y)$
can be expanded into power series by $y$, with the expansion
coefficients to be calculated by the following formulas:
$$
\gather
(a\circ S)_n=\Sum_{k=1}^n \Sum_{j_1+j_2+\dots+j_k=n}
a_k(S_{j_1},\dots,S_{j_k})\\
(b\circ S)_n=\Sum_{k=1}^n \Sum_{j_1+j_2+\dots+j_k=n}
b_k(S_{j_1},\dots,S_{j_k}).
\endgather
$$
Substituting expansions for $a(y)$ and $b(y)$ into original equation
and comparing corresponding coefficients, one can obtain the system of
linear stochastic equations for $n$-linear continuous operators
$S_n$, $n\geqslant1$:
$$ \cases S_1(t)&=S_1(0)+\Int_0^t a_1(\tau)S_1(\tau)d\tau+
\Int_0^t b_1(\tau)(S_1(\tau),dw(\tau)),\\ S_2(t)&=\Int_0^t
a_1(\tau)S_2(\tau)d\tau+ \Int_0^t b_1(\tau)(S_2(\tau),dw(\tau))+\\
&+\Int_0^t a_2(\tau)(S_1(\tau),S_1(\tau))d\tau+
\Int_0^t b_2(\tau)(S_1(\tau),S_1(\tau))dw(\tau),\\
&\qquad\vdots\\
S_n(t)&=\Int_0^t \Sum_{k=1}^n\Sum_{j_1+j_2+\dots+j_k=n}
a_k(\tau)(S_{j_1}(\tau),S_{j_2}(\tau),\dots,S_{j_k}(\tau))d\tau+\\
&+\Int_0^t \Sum_{k=1}^n\Sum_{j_1+j_2+\dots+j_k=n}
b_k(\tau)(S_{j_1}(\tau),S_{j_2}(\tau),\dots,S_{j_k}(\tau))dw(\tau),\\
&\qquad \vdots
\endcases
\tag1
$$
It's easy to see that the first $n$ equations of system \thetag1 with any
$n$ are closed with respect to $S_k$, $1\leqslant k\leqslant n$. It gives us a
possibility to solve the system recursively: find $S_1$ from the first
equation, find $S_2$ from the second one, using $S_1$ just found, find
$S_3$ from the third equations, using $S_1$ and $S_2$ already found, so on.

Let's note that system \thetag1 remains valid in the case, when $a_k$ and
$b_k$ ($k\geqslant 1$) are not coefficients of expansion of analytical function 
into
power series: all sums contained in \thetag1 are finite, and we may
implement recursion procedure without any worrying about series
convergence. These considerations lead to the notion of formal mapping.

\head 2. Formal mappings \endhead

Let $Y$ and $Z$ be Hilbert spaces.

\definition {Definition 1}
A sequence $a=(a_k)_{k\geqslant1}$, where $a_k$
($k\geqslant1$) are $k$-linear continuous mappings from $Y$ into $Z$,
we call formal mapping from $Y$ into $Z$.
Denote by $\CL_{\infty}(Y,Z)$ a space of formal mappings from $Y$
into $Z$.
\enddefinition

For formal mappings $a\in \CL_\infty(X,Y)$ and $b\in
\CL_\infty(Y,Z)$ composition operation is introduced:
$b\circ a\in\CL_\infty(X,Z)$,
$(b\circ a)_n=\Sum_{k=1}^n\Sum_{j_1+\dots+j_k=n} b_k
(a_{j_1},a_{j_2},\dots,a_{j_k})$.

Formal mapping $a\in \CL_\infty(Y,Y)$, denoted by $\Id_Y$, is called
identical one if $a_1=~\id_Y$ and $a_n=0$ with $n\geqslant2$.

\example {Example}
Any analytical function $a$ can be associated with formal mapping  $a$,
whose coefficients $a_k$ are coefficients of expanding $a$ into power series:
$a(y)=\break =\Sum_{k\geqslant1}a_k(y,y,\dots,y)$. In this case, composition of 
analytical
functions $c(y)=b(a(y))$ corresponds to composition of formal mappings
$c=b\circ a$, and identical function $a(y)=y$ corresponds to identical
formal mapping $\Id_Y$.
\endexample

To find more about formal mappings see \cite2.

\head
3. Stochastic equations in the space of formal mappings
\endhead

Let's consider stochastic equation in the space of formal mappings:
$$
S(t,s)=S(s,s)+\Int_s^t a(\tau)(S(\tau,s))d\tau
+\Int_s^t b(\tau)(S(\tau,s))dw(\tau),
\tag2
$$
where $a$ and $b$ are continuous mappings from
$[0,T]$ into spaces $\CL_\infty(Y,Y)$ and
$\CL_\infty(Y,\CL_2(Y,H_0))$ respectively, $S(t,s)$ is unknown random
process taking its values in space $\CL_\infty(Y,Y)$,
$S(s,s)$ is initial condition, measurable with respect to
$\sigma$-algebra $\F_s=\sigma(w(\tau),0\leqslant\tau\leqslant s)$.  Equation
\thetag2 is considered component-wise, i.e. \thetag2 is equivalent to
the following system:
$$
\cases
S_1(t,s)&=S_1(s,s)+\Int_s^t a_1(\tau)S_1(\tau,s)d\tau+ \Int_s^t
b_1(\tau)(S_1(\tau,s),dw(\tau)),\\
n\geqslant2:&\ S_n(t,s)=S_n(s,s)+\\
&+\Int_s^t \Sum_{k=1}^n\Sum_{j_1+j_2+\dots+j_k=n}
a_k(\tau)(S_{j_1}(\tau,s),S_{j_2}(\tau,s),\dots,S_{j_k}(\tau,s))d\tau+\\
&+\Int_s^t
\Sum_{k=1}^n\Sum_{j_1+j_2+\dots+j_k=n}
b_k(\tau)(S_{j_1}(\tau,s),S_{j_2}(\tau,s),\dots,S_{j_k}(\tau,s))dw(\tau).
\endcases
\tag3
$$

\proclaim {Theorem 1}
Let $a_n$ and $b_n$ ($n\geqslant1$) be functions, continuous with respect to 
$t\in
[0,T]$, taking their values in spaces $\CL_2(Y^{\otimes
n},Y)$ and $\CL_2(Y^{\otimes n}\otimes H_0,Y)$ respectively. Additionally
we suppose that the initial conditions satisfy the following
requirements:
$$ \gather S_1(s,s)-\id_Y\in\CL_2(Y,Y);\qquad
\forall\, n\geqslant2\: S_n(s,s)\in \CL_2(Y^{\otimes n},Y).
\endgather
$$
Then system \thetag2 has a solution $S(t,s)$, unique within stochastic
equivalence, such that:
$$ \gather
S_1(t,s)-\id_Y\in\CL_2(Y,Y); \qquad \forall\, n\geqslant2\: S_n(t,s)\in
\CL_2(Y^{\otimes n},Y).
\endgather
$$
\endproclaim

\proclaim{Theorem 2}
Let the condition of Theorem 1 hold, and operators
$S(t,s)$ ($0\leqslant s\leqslant t\leqslant T$) be solutions to equation
\thetag2 with initial condition $S(s,s)=\Id_Y$.

Then operator family $\{S(t,s)\}_{0\leqslant s\leqslant t\leqslant T}$ is 
evolution one,
i.e.:
$$
S(s,s)=\Id_Y,\quad S(t,\tau)\circ S(\tau,s)=S(t,s) \text{ with } s\leqslant
\tau \leqslant t.
$$
\endproclaim

Further we suppose that formal mapping $S(t)=S(t,0)$
is solution to \thetag2 with initial condition
$S(0)=\Id_Y$.
Let's rewrite system \thetag3 for process $S(t)$:
$$
\cases
S_1(t)&=\id_Y+\Int_0^t a_1(\tau)S_1(\tau)d\tau+ \Int_0^t
b_1(\tau)(S_1(\tau,s),dw(\tau)),\\
S_n(t)&=\Int_0^t a_1(\tau)S_n(\tau)d\tau+ \Int_0^t
b_1(\tau)(S_n(\tau),dw(\tau))+\\
&+\Int_0^t f_n(\tau)d\tau+\Int_0^t g_n(\tau)dw(\tau),\quad n\geqslant2,
\endcases
$$
\vskip-5pt
$$
\align
\text{where } f_n(t)&=\Sum_{k=2}^n\Sum_{j_1+j_2+\dots+j_k=n}
a_k(t) (S_{j_1}(t),S_{j_2}(t),\dots,S_{j_k}(t)),\\
g_n(t)&=\Sum_{k=2}^n\Sum_{j_1+j_2+\dots+j_k=n}
b_k(t) (S_{j_1}(t),S_{j_2}(t),\dots,S_{j_k}(t)).
\endalign
$$
Inasmuch as $f_n$ and $g_n$ contains only $S_k$ with $k\leqslant n-1$, this
system can be solved recursively, calculating $S_1$,
$S_2$, \dots, $S_n$, \dots. To give recursion procedure in a more
convenient way, we can use one explicit formula for solution to linear
nonhomogeneous stochastic equation (see \cite2):
$$ S_n(t)=\Int_0^t
S_1(t,\tau)(f_n(\tau))d\tau+ \Int_s^t
S_1(t,\tau)(g_n(\tau),\widehat{dw(\tau)}),
$$
where $S_1(t,\tau)$ is evolution operator, satisfying linear equation for
$S_1$, and symbol $\widehat{dw(\tau)}$ denotes extended stochastic
integral, treated as adjoint to stochastic derivative operator.

\bigskip
{\smc I.Ya. Spectorsky, postgraduate student, National Technical
University of Uk\-ra\-i\-ne "Kiev Polytechnic Institute", Chair of
Mathematical Methods of System Analysis, pr. Pobedy, 37, Kiev, Ukraine.
}

\enddocument